\documentclass[11pt,a4paper]{article}
\usepackage[T1]{fontenc}
\usepackage{lmodern}
\usepackage[margin=1in]{geometry}
\usepackage{amsmath,amssymb,amsthm}
\usepackage{enumitem}
\usepackage{microtype}
\usepackage[colorlinks=true,linkcolor=blue,citecolor=blue,urlcolor=blue]{hyperref}

\theoremstyle{plain}
\newtheorem{theorem}{Theorem}[section]
\newtheorem{lemma}[theorem]{Lemma}
\newtheorem{proposition}[theorem]{Proposition}
\newtheorem{corollary}[theorem]{Corollary}
\newtheorem{conjecture}[theorem]{Conjecture}
\theoremstyle{definition}
\newtheorem{definition}[theorem]{Definition}
\newtheorem{question}[theorem]{Question}
\theoremstyle{remark}
\newtheorem{remark}[theorem]{Remark}

\newcommand{\Rep}{\operatorname{Rep}}
\newcommand{\Hom}{\operatorname{Hom}}
\newcommand{\cC}{\mathcal C}
\newcommand{\cE}{\mathcal E}
\newcommand{\cG}{\mathcal G}
\newcommand{\Gop}{\mathcal G^{\mathrm{op}}}
\newcommand{\Ec}{\mathcal E_{\mathrm c}}

\title{The categorical graph minor conjecture and matching complexes}
\author{Jiangdong Ai\thanks{School of Mathematical Sciences and LPMC, Nankai University, Tianjin 300071, China. {\tt jd@nankai.edu.cn}. Partially supported by the National Natural Science Foundation of China (No. 12522117) and Fundamental and Interdisciplinary Disciplines Breakthrough Plan of the Ministry of Education of China (JYB2025XDXM207).}
\hspace{2mm}
    Yongtang Shi\thanks{School of Mathematical Sciences and LPMC, Nankai University, Tianjin 300071, China. {\tt shi@nankai.edu.cn}.
    Partially supported by the National Natural Science Foundation of China (No. 12431013).}
}
\date{}
\begin{document}
\maketitle

\begin{abstract}
In the graph minor category $\cG$ of Miyata, Proudfoot and Ramos, the morphisms are the ways of obtaining one connected graph from another by deleting and contracting edges. The categorical graph minor conjecture asserts that every submodule of a finitely generated $\Gop$-module is finitely generated; it is known for graphs of bounded first Betti number. We disprove this conjecture: it fails for the module assigning to each graph the free module on its spanning trees, which is generated by one element. A quotient of this module has unbounded torsion. We prove, however, that the conjecture holds for all modules depending on a minor morphism only through the induced injection of edge sets. This applies to matching complexes. Let $M_d(G)$ be the complex of edge sets of $G$ spanning subgraphs of maximum degree at most $d$. We prove that the $\Gop$-module $G\mapsto H_i(M_d(G);\mathbb Z)$ is finitely generated for all $i$ and $d$, which for $d=1$ confirms a conjecture recorded by Ramos, and that a positive integer depending only on $i$ annihilates the torsion of all these groups. The proofs combine Robertson and Seymour's labelled graph minor theorem with a new criterion: if a category has a faithful functor to finite sets and injections and, for every finite set of colours, its objects with coloured underlying sets are well-quasi-ordered by colour-preserving morphisms, then its representations over a left-Noetherian ring form a locally Noetherian category. Unlike the Gr\"obner method of Sam and Snowden, this criterion requires no ordering compatible with the morphisms.


\end{abstract}

\section{Introduction}\label{sec:intro}

The graph minor theorem of Robertson and Seymour~\cite{RS20} states that finite graphs are well-quasi-ordered by the minor relation: in every infinite sequence of finite graphs, some graph is isomorphic to a minor of a later one. Miyata, Proudfoot and Ramos~\cite{MPR} proposed a categorical form of this theorem, which we now describe. Throughout the paper, graphs are finite and may have loops and parallel edges.

The \emph{graph minor category} $\cG$ has the nonempty connected graphs as objects. A morphism from $H$ to $G$ in $\cG$, called a \emph{minor morphism}, is a way of obtaining $G$ from $H$ by contracting and deleting edges. It specifies which edges of $H$ are contracted, which are deleted, and how the remaining edges of $H$ are identified with the edges of $G$. The precise definition is recalled in Section~\ref{sec:trees}. Let $R$ be a ring. A \emph{$\Gop$-module} over $R$ is a functor $M$ from the opposite category $\Gop$ to the category of left $R$-modules. Thus every minor morphism $\phi$ from $H$ to $G$ induces a homomorphism $M(\phi)\colon M(G)\to M(H)$. A \emph{submodule} $N$ of $M$ consists of submodules $N(G)\subseteq M(G)$ which are preserved by these homomorphisms. The module $M$ is \emph{finitely generated} if it has finitely many elements $x_j\in M(G_j)$ that are contained in no proper submodule of $M$.

The main result announced in~\cite{MPR} was that every submodule of a finitely generated $\Gop$-module over a Noetherian ring is finitely generated. A gap was later found in the proof~\cite[Remark~1.5]{MR}, and the statement has since been treated as a conjecture, the categorical graph minor conjecture.

\begin{conjecture}[{\cite[Conjecture~2.3]{RamosSurvey}, \cite[Conjecture~4.1]{KR}}]\label{conj:cgm}
Let $R$ be a Noetherian ring. Then every submodule of a finitely generated $\Gop$-module over $R$ is finitely generated.
\end{conjecture}

Miyata and Ramos~\cite[Conjecture~3.3]{MR} conjectured the stronger statement that $\Gop$ is quasi-Gr\"obner in the sense of Sam and Snowden~\cite{SS}. Conjecture~\ref{conj:cgm} is known for graphs of bounded first Betti number. Recall that the first Betti number of a connected graph $G$ is $|E(G)|-|V(G)|+1$; it is called the genus of $G$ in~\cite{MR}. For an integer $g\ge0$, let $\cG_{\le g}$ be the full subcategory of $\cG$ on the graphs with first Betti number at most $g$. Building on the work of Proudfoot and Ramos~\cite{PR}, Miyata and Ramos proved that $\cG_{\le g}^{\mathrm{op}}$ is quasi-Gr\"obner~\cite[Theorem~3.5]{MR}, and hence that every submodule of a finitely generated $\cG_{\le g}^{\mathrm{op}}$-module over a Noetherian ring is finitely generated~\cite[Theorem~1.1]{MR}. Analogous results for quivers were obtained by Caputi, Collari and Ramos~\cite{CCR}.

Conjecture~\ref{conj:cgm} was intended as a tool for proving finiteness theorems which are uniform over all graphs, for instance for the homology of matching complexes~\cite{MR,RamosSurvey} and of configuration spaces of graphs~\cite{KR}. In this paper, we show that Conjecture~\ref{conj:cgm} is false, we prove a weaker statement of the same kind, and we show that the weaker statement suffices for matching complexes.

\subsection*{A counterexample}

Our first result is the following.

\begin{theorem}\label{thm:trees}
Let $R$ be a nonzero ring. Then there is a $\Gop$-module over $R$ which is generated by a single element and has a submodule that is not finitely generated. Consequently, Conjecture~\ref{conj:cgm} is false, and $\Gop$ is not quasi-Gr\"obner.
\end{theorem}

The module in Theorem~\ref{thm:trees} assigns to a connected graph $G$ the free $R$-module on the set of spanning trees of $G$. A minor morphism from $H$ to $G$ maps a spanning tree $T$ of $G$ to the spanning tree of $H$ formed by the edges of $H$ which are identified with the edges of $T$, together with the contracted edges of $H$. This is the spanning tree module of~\cite[Section~2]{RamosSurvey}; it is the principal projective module of the graph with one vertex and no edges~\cite[Remark~3.19]{MR}. We show that the submodule generated by the spanning stars of the wheels is not finitely generated. All graphs in this example are planar. Their first Betti numbers are unbounded, in accordance with~\cite[Theorem~1.1]{MR}. The wheels can be replaced by graphs with no minor isomorphic to $K_4$ (Remark~\ref{rem:brooms}).

A quotient of the spanning tree module shows that the torsion of a finitely generated $\Gop$-module need not be bounded. In Proposition~\ref{prop:torsion} we construct a $\Gop$-module over $\mathbb Z$ which is generated by a single element and whose value at the wheel with $m$ spokes has a direct summand isomorphic to $\mathbb Z/m\mathbb Z$, for every $m\ge3$. A uniform bound on the torsion of finitely generated $\Gop$-modules is stated in~\cite[Section~2]{RamosSurvey} with reference to~\cite{MPR}, and bounds of this kind were among the main applications of the statement of Conjecture~\ref{conj:cgm}.

\subsection*{A categorical graph minor theorem for edge maps}

Every minor morphism $\phi$ from $H$ to $G$ induces an injection $\phi^*\colon E(G)\to E(H)$, which maps every edge of $G$ to the edge of $H$ identified with it. The module of spanning trees depends on more than this injection. For instance, all minor morphisms from $H$ to the graph with one vertex and no edges induce the empty injection, but they correspond to the different spanning trees of $H$. Our second result shows that Conjecture~\ref{conj:cgm} holds for the modules which depend on a minor morphism $\phi$ only through $\phi^*$.

\begin{theorem}\label{thm:edge}
Let $R$ be a left-Noetherian ring, and let $M$ be a finitely generated $\Gop$-module over $R$. Suppose that $M(\phi)=M(\phi')$ whenever $\phi$ and $\phi'$ are minor morphisms from $H$ to $G$ with $\phi^*=\phi'^*$. Then every submodule of $M$ is finitely generated.
\end{theorem}

Theorem~\ref{thm:edge} can be restated as follows. The injections $\phi^*$ are the morphisms of a category $\Ec$ on the nonempty connected graphs, and a $\Gop$-module as in Theorem~\ref{thm:edge} is the same as a functor from $\Ec$ to the category of left $R$-modules. So Theorem~\ref{thm:edge} says that every subfunctor of a finitely generated functor of this kind is finitely generated. We prove this for the larger category $\cE$ of all graphs, whose morphisms from $G$ to $H$ are the injections $E(G)\to E(H)$ induced by minor models of $G$ in $H$ (Proposition~\ref{prop:E}). Like Conjecture~\ref{conj:cgm}, this statement is a categorical form of the graph minor theorem, which it implies (Remark~\ref{rem:gmt}). Its proof relies on the labelled graph minor theorem of Robertson and Seymour~\cite[1.7]{RS23}, and so it does not give a new proof of the graph minor theorem.

\subsection*{Matching complexes}

Let $G$ be a graph and let $d\ge1$ be an integer. For a set $F\subseteq E(G)$ and a vertex $v$ of $G$, the \emph{degree} of $v$ in $F$ is the number of ends of edges of $F$ at $v$, so that a loop at $v$ contributes two. The \emph{$d$-matching complex} $M_d(G)$ is the simplicial complex on $E(G)$ whose faces are the sets $F\subseteq E(G)$ in which every vertex of $G$ has degree at most $d$. In particular, the faces of $M_1(G)$ are the matchings of $G$, that is, the sets of non-loop edges no two of which share an end, and $M_1(G)$ is the matching complex of $G$, as in~\cite[Section~4.1]{MR}. The homology of matching complexes has a rich structure; already for complete graphs it contains torsion, as discovered by Bouc~\cite{Bouc}. We refer to Shareshian and Wachs~\cite{SW} for these results.

For a minor morphism $\phi$ from $H$ to $G$, the injection $\phi^*$ maps the faces of $M_d(G)$ to faces of $M_d(H)$. Hence $G\mapsto H_i(M_d(G);R)$ is a $\Gop$-module for all integers $i\ge0$ and $d\ge1$, where $H_i(\,\cdot\,;R)$ denotes simplicial homology with coefficients in $R$, and this module depends on a minor morphism $\phi$ only through $\phi^*$. Miyata and Ramos proved the following.

\begin{theorem}[{\cite[Theorem~4.3 and Corollary~4.4]{MR}}]\label{thm:MR}
Let $g,i\ge0$ and $d\ge1$ be integers. Then the $\cG_{\le g}^{\mathrm{op}}$-module $G\mapsto H_i(M_d(G);\mathbb Z)$ is finitely generated. Moreover, there is a positive integer $\varepsilon_{i,d,g}$ which annihilates the torsion subgroup of $H_i(M_d(G);\mathbb Z)$ for every connected graph $G$ with first Betti number at most $g$.
\end{theorem}

In his survey, Ramos records the conjecture that, for matching complexes, the first assertion of Theorem~\ref{thm:MR} holds without the restriction on the first Betti number.

\begin{conjecture}[{\cite[Section~2]{RamosSurvey}}]\label{conj:matching}
For every integer $i\ge0$, the $\Gop$-module $G\mapsto H_i(M_1(G);\mathbb Z)$ is finitely generated.
\end{conjecture}

Conjecture~\ref{conj:matching} was expected to follow from Conjecture~\ref{conj:cgm}, by the proof of Theorem~\ref{thm:MR}; see~\cite[Section~3]{RamosSurvey} and~\cite[Remark~3.4]{MR}. By Theorem~\ref{thm:trees} this route is not available. We prove Conjecture~\ref{conj:matching} by the method of Theorem~\ref{thm:edge}, for all $d$ and together with a bound on the torsion.

\begin{theorem}\label{thm:matching}
Let $i\ge0$ be an integer.
\begin{enumerate}[label=\textup{(\roman*)}]
\item Let $d\ge1$ be an integer and let $R$ be a left-Noetherian ring. Then the $\Gop$-module $G\mapsto H_i(M_d(G);R)$ is finitely generated.
\item There is a positive integer $b_i$ such that $b_i$ annihilates the torsion subgroup of $H_i(M_d(G);\mathbb Z)$ for every graph $G$ and every integer $d\ge1$.
\end{enumerate}
\end{theorem}

Part~(i) with $d=1$ and $R=\mathbb Z$ is Conjecture~\ref{conj:matching}. It says that there are finitely many connected graphs $G_1,\dots,G_n$ and homology classes $x_j\in H_i(M_1(G_j);\mathbb Z)$ such that, for every connected graph $H$, the group $H_i(M_1(H);\mathbb Z)$ is generated by the images of the classes $x_j$ under the homomorphisms induced by the minor morphisms from $H$ to the graphs $G_j$. In part~(ii) the graph $G$ need not be connected, and $b_i$ depends neither on the first Betti number of $G$ nor on $d$. Any such integer $b_1$ is divisible by $3$, since $H_1(M_1(K_7);\mathbb Z)\cong\mathbb Z/3\mathbb Z$ by a result of Bouc~\cite{Bouc}; see~\cite[Theorem~1.3(ii)]{SW}. Our proof gives no explicit bound on $b_i$. The same method applies to other complexes of sets of edges (Corollary~\ref{cor:properties}).

\subsection*{A criterion for Noetherianity}

Sam and Snowden~\cite{SS} showed how Noetherian properties of representations of combinatorial categories can be derived from well-quasi-orders. Their Gr\"obner method requires, besides a well-quasi-order on the morphisms out of each object, an ordering of these morphisms which is compatible with composition, possibly after passing to a suitable cover of the category. For graphs of bounded first Betti number such orderings are available~\cite{MR,PR}. By Theorem~\ref{thm:trees} there is no suitable ordered cover of $\Gop$, and we do not know of one for the category $\cE$. Our proof of Theorem~\ref{thm:edge} is based on a criterion in which no compatible ordering is required. Instead, the well-quasi-order is assumed for all colourings with finitely many colours.

Let $\cC$ be an essentially small category equipped with a faithful functor $V$ from $\cC$ to the category of finite sets and injections. A \emph{representation} of $\cC$ over $R$ is a functor from $\cC$ to the category of left $R$-modules, and $\Rep_R(\cC)$ denotes the category of these representations. Subrepresentations and finitely generated representations are defined as for $\Gop$-modules above. We say that $\Rep_R(\cC)$ is \emph{locally Noetherian} if every subrepresentation of a finitely generated representation is finitely generated. For a finite set $D$, a \emph{$D$-coloured object} of $\cC$ is a pair $(y,c)$ consisting of an object $y$ and a map $c\colon V(y)\to D$. A morphism $h\colon y\to z$ of $\cC$ is \emph{colour-preserving} from $(y,c)$ to $(z,c')$ if $c'(V(h)(v))=c(v)$ for every $v\in V(y)$. A quasi-order $(Q,\le)$ is a \emph{well-quasi-order}, or \emph{wqo}, if every infinite sequence $q_1,q_2,\dots$ of elements of $Q$ contains indices $j<k$ with $q_j\le q_k$. Accordingly, the $D$-coloured objects of $\cC$ are well-quasi-ordered by colour-preserving morphisms if every infinite sequence of $D$-coloured objects contains indices $j<k$ such that there is a colour-preserving morphism from the $j$th term to the $k$th term.

\begin{theorem}\label{thm:criterion}
Suppose that, for every finite set $D$, the $D$-coloured objects of $\cC$ are well-quasi-ordered by colour-preserving morphisms. Then $\Rep_R(\cC)$ is locally Noetherian for every left-Noetherian ring $R$. The same holds for every full subcategory of $\cC$.
\end{theorem}

The proof is a leading-term argument. We order the morphisms from a fixed object $x$ to each object $y$ by means of an arbitrary linear order on $V(y)$, chosen independently for every $y$. Given a morphism $f\colon x\to y$, where $V(x)$ has $m$ elements, we give the $m$ elements of the image of $f$ distinct colours and colour every other element of $V(y)$ by the number of elements of the image that precede it. A morphism preserving these $2m+1$ colours need not preserve the linear orders, but it preserves all comparisons with $f$, and this is what the elimination of a leading term requires. The argument is related to the nice enumerations of Ahlbrandt and Ziegler~\cite{AZ}, which Camina and Evans~\cite{CE} used to prove the ascending chain condition for submodules of permutation modules, and to the symmetric cancellation ordering of Aschenbrenner and Hillar~\cite{AH}. Ghosh and Lasota~\cite[Theorem~27]{GL} proved that chains of equivariant subspaces of tuple spaces are finite under an ordered expansion hypothesis, and Laudone and Snowden~\cite{LS} studied Noetherianity for systems of polynomial ideals on combinatorial categories using compatible orders, possibly on a cover.

For the category $\cE$ and the functor $G\mapsto E(G)$, the hypothesis of Theorem~\ref{thm:criterion} is the edge-labelled graph minor theorem of Robertson and Seymour~\cite[1.7]{RS23}. For connected graphs every morphism of $\cE$ is of the form $\phi^*$ (Lemma~\ref{lem:full}), and this gives Theorem~\ref{thm:edge}. Theorem~\ref{thm:criterion} applies to other well-quasi-orders of graphs as well. With vertex sets in place of edge sets, the labelled immersion theorem of Robertson and Seymour~\cite[1.5]{RS23} shows that the categories of vertex-labelled graphs whose morphisms are the vertex maps of immersions have locally Noetherian representations (Corollary~\ref{cor:immersion}). For comparison, Knudsen and Ramos~\cite[Theorem~7.2]{KR} proved that the category of cographs and full embeddings is Noetherian, and Mamun, Nalikka and Ramos~\cite{MNR} obtained related results for polynomial-ring modules on cographs.

\subsection*{Organisation}

The paper is organised as follows. In Section~\ref{sec:trees} we recall the graph minor category and prove Theorem~\ref{thm:trees}. Section~\ref{sec:criterion} contains the proof of Theorem~\ref{thm:criterion} and a corollary for simplicial complexes. In Section~\ref{sec:minors} we study the category $\cE$ and prove Theorem~\ref{thm:edge}. Theorem~\ref{thm:matching} is proved in Section~\ref{sec:matching}. Section~\ref{sec:immersions} deals with immersions, and Section~\ref{sec:conclusion} contains concluding remarks and open questions.

\section{The graph minor category and spanning trees}\label{sec:trees}

In this section we recall the definition of the graph minor category and prove Theorem~\ref{thm:trees}.

\subsection{The graph minor category}

Recall that graphs are finite and may have loops and parallel edges. Every edge has two ends, which coincide if the edge is a loop. A subgraph is \emph{connected} if it is nonempty and any two of its vertices are joined by a path in it. A \emph{tree} is a connected graph without cycles; in particular, a tree has no loops and no parallel edges.

The following definition is taken from~\cite[Section~3.1]{MR}. Every edge of a graph is regarded as a pair of opposite \emph{arrows}. An arrow has a head and a tail, which are the two ends of its edge, and reversing an arrow exchanges its head and its tail. We write $A(G)$ for the set of arrows of a graph $G$, and we fix an additional symbol $\star$. Let $G$ and $H$ be graphs. A \emph{minor morphism} from $H$ to $G$ is a map $\phi$ from $V(H)\cup A(H)\cup\{\star\}$ to $V(G)\cup A(G)\cup\{\star\}$ with the following properties.
\begin{enumerate}[label=\textup{(M\arabic*)}]
\item $\phi(\star)=\star$, and $\phi$ maps vertices to vertices.
\item $\phi$ commutes with the reversal of arrows, where vertices and $\star$ are regarded as their own reversals.
\item Every arrow of $G$ is the image of exactly one arrow of $H$.
\item If $\phi$ maps an arrow $a$ to an arrow, then it maps the head and the tail of $a$ to the head and the tail of $\phi(a)$.
\item If $\phi$ maps an arrow $a$ to a vertex, then it maps the head and the tail of $a$ to this vertex.
\item For every vertex $v$ of $G$, the vertices and edges of $H$ which are mapped to $v$ form a tree.
\end{enumerate}
By~(M2), the two arrows of an edge of $H$ are both mapped to the same vertex, or both to $\star$, or to the two arrows of an edge of $G$. In the first case the edge is \emph{contracted} by $\phi$, in the second case it is \emph{deleted}, and in the third case it is \emph{mapped to an edge}. We write $K(\phi)$ for the set of edges of $H$ contracted by $\phi$. Minor morphisms are composed as maps. The \emph{graph minor category} $\cG$ has the nonempty connected graphs as objects and the minor morphisms as morphisms.

Let $\phi$ be a minor morphism from $H$ to $G$. By~(M2) and~(M3), every edge $e$ of $G$ is the image of exactly one edge of $H$, which we denote by $\phi^*(e)$. This defines an injection $\phi^*\colon E(G)\to E(H)$.

\subsection{The module of spanning trees}

Let $\bullet$ denote the graph with one vertex and no edges, and let $G$ be a nonempty connected graph. A minor morphism $\tau$ from $G$ to $\bullet$ maps every vertex of $G$ to the vertex of $\bullet$, and it contracts or deletes every edge of $G$. By~(M6), the vertices of $G$ and the contracted edges form a tree, that is, $K(\tau)$ is the edge set of a spanning tree of $G$. Conversely, for every spanning tree $T$ of $G$ there is exactly one minor morphism from $G$ to $\bullet$ which contracts the edges of $T$ and deletes all other edges. So $\tau\mapsto K(\tau)$ is a bijection between $\Hom_{\cG}(G,\bullet)$ and the set of spanning trees of $G$; compare~\cite[Remark~3.19]{MR}.

Let $R$ be a ring. Let $P(G)$ be the free left $R$-module on the set $\Hom_{\cG}(G,\bullet)$, and for a minor morphism $\phi$ from $H$ to $G$ let $P(\phi)\colon P(G)\to P(H)$ be the homomorphism mapping a basis element $\tau$ to $\tau\phi$. This defines a $\Gop$-module $P$ over $R$, the principal projective $\Gop$-module of the object $\bullet$. An edge of $H$ is contracted by $\tau\phi$ if and only if it is contracted by $\phi$, or mapped by $\phi$ to an edge that is contracted by $\tau$. Hence
\[
K(\tau\phi)=K(\phi)\cup\phi^*(K(\tau)).
\]
In terms of spanning trees, $P(G)$ is the free module on the spanning trees of $G$, and $P(\phi)$ maps a spanning tree $T$ of $G$ to the spanning tree of $H$ with edge set $K(\phi)\cup\phi^*(E(T))$. This is the module described in the introduction.

For $k\ge3$ let $W_k$ be the wheel with $k$ spokes. It has a vertex $z$, called the centre, and $k$ further vertices $u_1,\dots,u_k$. Its edges are the \emph{spokes} $zu_j$ and the \emph{rim edges} $u_ju_{j+1}$ for $1\le j\le k$, where $u_{k+1}=u_1$. The rim edges form a cycle of length $k$, called the rim. Let $\sigma_k$ be the minor morphism from $W_k$ to $\bullet$ which contracts the spokes and deletes the rim edges.

\begin{lemma}\label{lem:wheels}
Let $3\le k<m$. There is no minor morphism $\psi$ from $W_m$ to $W_k$ with $\sigma_k\psi=\sigma_m$.
\end{lemma}

\begin{proof}
Suppose that $\psi$ is such a minor morphism. Let $S_k$ and $S_m$ denote the sets of spokes of $W_k$ and $W_m$. Then
\[
S_m=K(\sigma_m)=K(\sigma_k\psi)=K(\psi)\cup\psi^*(S_k).
\]
So every edge contracted by $\psi$ is a spoke of $W_m$, and $\psi^*(e)$ is a spoke of $W_m$ for every spoke $e$ of $W_k$. If $e$ is a rim edge of $W_k$, then $\psi^*(e)$ is not contracted by $\psi$ and does not belong to $\psi^*(S_k)$, because $\psi^*$ is injective. So $\psi^*(e)$ is not a spoke, that is, $\psi^*(e)$ is a rim edge of $W_m$.

For a vertex $w$ of $W_k$, the vertices and edges of $W_m$ mapped to $w$ form a tree by~(M6), and all edges of this tree are spokes of $W_m$. Every spoke of $W_m$ is incident with the centre $z_m$ of $W_m$, so such a tree consists of a single vertex unless it contains $z_m$. Let $w_0=\psi(z_m)$. Each of the $k$ vertices of $W_k$ other than $w_0$ is therefore the image of exactly one vertex of $W_m$, and this vertex is not $z_m$. Let $U$ be the set of these $k$ vertices of $W_m$. They lie on the rim of $W_m$.

Let $e$ be a spoke of $W_k$. Then $\psi^*(e)$ is a spoke of $W_m$, so $z_m$ is one of its ends, and by~(M4) the vertex $w_0=\psi(z_m)$ is an end of $e$. Thus $w_0$ is incident with all $k\ge3$ spokes of $W_k$, and so $w_0$ is the centre of $W_k$. Now let $e$ be a rim edge of $W_k$. Its ends are distinct from $w_0$, and $\psi^*(e)$ is a rim edge of $W_m$ whose ends are mapped to the ends of $e$ by~(M4). So both ends of $\psi^*(e)$ belong to $U$. As $\psi^*$ is injective, the rim of $W_m$ has at least $k$ edges with both ends in $U$. But the rim of $W_m$ is a cycle of length $m>k$, and a set of $k$ of its vertices spans at most $k-1$ of its edges. This is a contradiction.
\end{proof}

\begin{proof}[Proof of Theorem~\ref{thm:trees}]
Let $R$ be a nonzero ring and let $P$ be the $\Gop$-module defined above. For every minor morphism $\tau$ from $G$ to $\bullet$, the homomorphism $P(\tau)$ maps the identity of $\bullet$ to $\tau$. So $P$ is generated by a single element. Let $N$ be the submodule of $P$ generated by the elements $\sigma_k\in P(W_k)$ for $k\ge3$. For every graph $G$, the module $N(G)$ is spanned by the elements $\sigma_k\psi$, where $k\ge3$ and $\psi$ is a minor morphism from $G$ to $W_k$.

Suppose that $N$ is finitely generated. Each of finitely many generators is a linear combination of finitely many elements $\sigma_k\psi$, so there is an integer $m>3$ such that $N$ is generated by the elements $\sigma_k$ with $3\le k<m$. Then $\sigma_m$ is a linear combination of elements $\sigma_k\psi$, where $3\le k<m$ and $\psi$ is a minor morphism from $W_m$ to $W_k$. Both $\sigma_m$ and the elements $\sigma_k\psi$ belong to the basis $\Hom_{\cG}(W_m,\bullet)$ of the free module $P(W_m)$. As $R$ is nonzero, it follows that $\sigma_m=\sigma_k\psi$ for some $k$ and $\psi$, which contradicts Lemma~\ref{lem:wheels}. Hence $N$ is not finitely generated.

Taking $R$ to be a nonzero Noetherian ring, we see that Conjecture~\ref{conj:cgm} is false. If $\Gop$ were quasi-Gr\"obner, then every submodule of a finitely generated $\Gop$-module over a left-Noetherian ring would be finitely generated by~\cite[Theorem~1.1.3]{SS}. So $\Gop$ is not quasi-Gr\"obner.
\end{proof}

\begin{proposition}\label{prop:torsion}
There is a $\Gop$-module $M$ over $\mathbb Z$ which is generated by a single element and such that, for every $m\ge3$, the group $M(W_m)$ has a direct summand isomorphic to $\mathbb Z/m\mathbb Z$. In particular, no positive integer annihilates the torsion subgroups of all the groups $M(G)$.
\end{proposition}

\begin{proof}
Let $R=\mathbb Z$, and let $N'$ be the submodule of $P$ generated by the elements $k\sigma_k\in P(W_k)$ for $k\ge3$. Let $M$ be the quotient of $P$ by $N'$, so that $M(G)=P(G)/N'(G)$ for every nonempty connected graph $G$. As $P$ is generated by the identity of $\bullet$, the module $M$ is generated by the image of this element.

Let $m\ge3$. The group $N'(W_m)$ is generated by the elements $k\,\sigma_k\psi$, where $k\ge3$ and $\psi$ is a minor morphism from $W_m$ to $W_k$, and each element $\sigma_k\psi$ belongs to the basis $\Hom_{\cG}(W_m,\bullet)$ of $P(W_m)$. For $\tau\in\Hom_{\cG}(W_m,\bullet)$ let $n_\tau$ be the greatest common divisor of the integers $k\ge3$ such that $\tau=\sigma_k\psi$ for some minor morphism $\psi$ from $W_m$ to $W_k$, where $n_\tau=0$ if there is no such $k$. Then $N'(W_m)$ is the direct sum of the subgroups $n_\tau\mathbb Z\,\tau$ of $P(W_m)$, and so $M(W_m)$ is isomorphic to the direct sum of the groups $\mathbb Z/n_\tau\mathbb Z$. It remains to show that $n_\tau=m$ for $\tau=\sigma_m$. The integer $k=m$ occurs, with the identity of $W_m$ as $\psi$. No integer $k$ with $3\le k<m$ occurs, by Lemma~\ref{lem:wheels}. If $k>m$, then there is no minor morphism from $W_m$ to $W_k$ at all, because by~(M6) every vertex of $W_k$ would be the image of a vertex of $W_m$, and $W_k$ has more vertices than $W_m$.
\end{proof}

\begin{remark}\label{rem:wheels}
The proof of Theorem~\ref{thm:trees} uses only the wheels and the graph $\bullet$. Hence the conclusion of Theorem~\ref{thm:trees} holds for every full subcategory of $\cG$ which contains $\bullet$ and infinitely many wheels, for instance for the full subcategory on the connected planar graphs. The module $P$ is not of the kind considered in Theorem~\ref{thm:edge}: for every connected graph $G$, all minor morphisms $\tau$ from $G$ to $\bullet$ induce the same injection, namely the empty one, whereas $P(\tau)$ maps the identity of $\bullet$ to $\tau$.
\end{remark}

\begin{remark}\label{rem:brooms}
The wheels in the proof of Theorem~\ref{thm:trees} can be replaced by graphs with no minor isomorphic to $K_4$. For $n\ge2$ let $F_n$ be the tree obtained from a path with $n$ vertices by attaching two new leaves to each of the two ends of the path, and let $A_n$ be the graph obtained from $F_n$ by adding a new vertex $z_n$ and one edge from $z_n$ to each vertex of $F_n$. These new edges are the spokes of $A_n$. Let $\rho_n$ be the minor morphism from $A_n$ to $\bullet$ which contracts the spokes and deletes the edges of $F_n$.

Let $2\le k<m$, and suppose that $\psi$ is a minor morphism from $A_m$ to $A_k$ with $\rho_k\psi=\rho_m$. The proof of Lemma~\ref{lem:wheels}, with the edges of $F_k$ and $F_m$ in place of the rim edges, shows that $\psi(z_m)=z_k$, that every other vertex of $A_k$ is the image of exactly one vertex of $A_m$, and that $\psi^*$ maps every edge $e$ of $F_k$ to an edge of $F_m$ whose ends are mapped to the ends of $e$. So $F_k$ is isomorphic to a subgraph of $F_m$. This is impossible: the two vertices of degree three of $F_k$ are joined by a path of length $k-1$, they would correspond to the only two vertices of degree at least three of $F_m$, and in the tree $F_m$ the only path joining these two vertices has length $m-1$. As in the proof of Theorem~\ref{thm:trees}, it follows that the submodule of $P$ generated by the elements $\rho_n$ with $n\ge2$ is not finitely generated.

The graph $A_n$ has no minor isomorphic to $K_4$: otherwise a triangle, which is what remains of $K_4$ after deleting a vertex, would be a minor of the tree $A_n-z_n=F_n$. Hence the conclusion of Theorem~\ref{thm:trees} also holds for the full subcategory of $\cG$ on the connected graphs with no minor isomorphic to $K_4$.
\end{remark}

\section{The finite-colouring criterion}\label{sec:criterion}

Throughout this section, $\cC$ is an essentially small category equipped with a faithful functor $V$ to the category of finite sets and injections, and $R$ is a ring. For a morphism $h\colon y\to z$ of $\cC$ and an element $v\in V(y)$ we write $h(v)$ for $V(h)(v)$. Since $V$ is faithful, two morphisms $h,h'\colon y\to z$ with $h(v)=h'(v)$ for every $v\in V(y)$ are equal.

\subsection{Preliminaries}

Let $M$ be a representation of $\cC$ over $R$. For a morphism $h\colon y\to z$ we write $M(h)\colon M(y)\to M(z)$ for the induced homomorphism. A \emph{subrepresentation} $N$ of $M$ consists of submodules $N(y)\subseteq M(y)$ with $M(h)(N(y))\subseteq N(z)$ for every morphism $h\colon y\to z$. The subrepresentation \emph{generated} by a family of elements $x_j\in M(y_j)$ is the smallest subrepresentation containing them. Its value at an object $z$ consists of the finite sums of elements $r\,M(h)(x_j)$ with $r\in R$ and $h\colon y_j\to z$. The representation $M$ is \emph{finitely generated} if it is generated by a finite family of elements.

For an object $x$, the \emph{principal projective} representation $P_x$ is defined as follows. The module $P_x(y)$ is the free left $R$-module on the set $\Hom_{\cC}(x,y)$, and for a morphism $h\colon y\to z$ the homomorphism $P_x(h)$ maps a basis element $g$ to $hg$. For every representation $M$ and every element $w\in M(x)$ there is a unique morphism of representations $P_x\to M$ which maps the identity of $x$ to $w$, namely the one mapping $g\in\Hom_{\cC}(x,y)$ to $M(g)(w)$. Its image is the subrepresentation generated by $w$.

We use the following standard facts on well-quasi-orders.

\begin{lemma}\label{lem:wqo}
Let $(Q,\le)$ be a wqo.
\begin{enumerate}[label=\textup{(\roman*)}]
\item Every infinite sequence $q_1,q_2,\dots$ of elements of $Q$ has an infinite subsequence $q_{j_1}\le q_{j_2}\le\cdots$ with $j_1<j_2<\cdots$.
\item Every subset $S$ of $Q$ has a finite subset $S_0$ such that every element of $S$ is greater than or equal to some element of $S_0$.
\end{enumerate}
\end{lemma}

\begin{proof}
(i) Call an index $j$ \emph{terminal} if there is no $k>j$ with $q_j\le q_k$. If there were infinitely many terminal indices $j_1<j_2<\cdots$, then the sequence $q_{j_1},q_{j_2},\dots$ would contradict the definition of a wqo. So there is an index $j_1$ which is larger than all terminal indices. Since $j_1$ is not terminal, there is $j_2>j_1$ with $q_{j_1}\le q_{j_2}$. Since $j_2$ is not terminal, there is $j_3>j_2$ with $q_{j_2}\le q_{j_3}$. Continuing in this way we obtain the required subsequence.

(ii) Suppose that there is no such set $S_0$. Then we can choose elements $s_1,s_2,\dots$ of $S$ one after another such that $s_k$ is not greater than or equal to any of $s_1,\dots,s_{k-1}$. This sequence contradicts the definition of a wqo.
\end{proof}

\subsection{An order on morphisms}\label{subsec:order}

A skeleton of $\cC$ is a small full subcategory which contains exactly one object from every isomorphism class. The hypothesis of Theorem~\ref{thm:criterion} passes from $\cC$ to each of its full subcategories, since a full subcategory has fewer coloured objects and the same colour-preserving morphisms between them. Moreover, restriction to a skeleton is an equivalence between the categories of representations, which preserves subrepresentations and finite generation. Hence we may, and do, assume in the proof of Theorem~\ref{thm:criterion} that $\cC$ is small.

For every object $y$ we fix a linear order $<_y$ on $V(y)$. We also fix an object $x$ and an enumeration $u_1,\dots,u_m$ of $V(x)$, where $m=|V(x)|$. Let $y$ be an object. For distinct morphisms $f,g\colon x\to y$, we write $g<f$ if the tuple $(g(u_1),\dots,g(u_m))$ precedes the tuple $(f(u_1),\dots,f(u_m))$ in the lexicographic order defined by $<_y$. As $V$ is faithful, distinct morphisms have distinct tuples, so this is a linear order on the finite set $\Hom_{\cC}(x,y)$. Note that only morphisms with the same target are compared.

For morphisms $f\colon x\to y$ and $f'\colon x\to z$, we write $f\unlhd f'$ if there is a morphism $h\colon y\to z$ such that $hf=f'$ and $hg<f'$ for every morphism $g\colon x\to y$ with $g<f$. Such a morphism $h$ is called a \emph{witness} for $f\unlhd f'$. The identity of $y$ is a witness for $f\unlhd f$, and if $h$ is a witness for $f\unlhd f'$ and $h'$ is a witness for $f'\unlhd f''$, then $h'h$ is a witness for $f\unlhd f''$. So $\unlhd$ is a quasi-order on the set of morphisms with source $x$.

\begin{lemma}\label{lem:order}
Suppose that $\cC$ satisfies the hypothesis of Theorem~\ref{thm:criterion}. Then $\unlhd$ is a well-quasi-order on the set of morphisms of $\cC$ with source $x$.
\end{lemma}

\begin{proof}
Let $f_1,f_2,\dots$ be an infinite sequence of morphisms with source $x$, where $f_j\colon x\to y_j$. For each $j$, the order $<_{y_j}$ arranges the $m$ distinct elements $f_j(u_1),\dots,f_j(u_m)$ in one of $m!$ possible ways. Passing to an infinite subsequence, we may assume that this arrangement is the same for all $j$; that is, for all $s,t\in\{1,\dots,m\}$ and all $j,k$ we have $f_j(u_s)<_{y_j}f_j(u_t)$ if and only if $f_k(u_s)<_{y_k}f_k(u_t)$.

Let $D$ be a set of $2m+1$ colours $r_1,\dots,r_m,d_0,\dots,d_m$. We define a colouring $c_j\colon V(y_j)\to D$ as follows. For $1\le s\le m$ let $c_j(f_j(u_s))=r_s$. For an element $v$ of $V(y_j)$ outside the image of $f_j$ let $c_j(v)=d_b$, where $b$ is the number of indices $s$ with $f_j(u_s)<_{y_j}v$. By the hypothesis, there are indices $j<k$ and a colour-preserving morphism $h$ from $(y_j,c_j)$ to $(y_k,c_k)$. We show that $h$ is a witness for $f_j\unlhd f_k$, which proves the lemma.

First, for each $s$ the element $h(f_j(u_s))$ has colour $r_s$, and $f_k(u_s)$ is the only element of $V(y_k)$ with this colour. Hence $hf_j(u_s)=f_k(u_s)$ for every $s$, and so $hf_j=f_k$.

Secondly, we claim that $h(v)<_{y_k}f_k(u_s)$ for every $v\in V(y_j)$ and every $s$ such that $v<_{y_j}f_j(u_s)$. If $v=f_j(u_t)$ for some $t$, then $h(v)=f_k(u_t)$, and the claim follows from the choice of the subsequence. Otherwise $v$ has colour $d_b$ for some $b$. Exactly $b$ of the elements $f_j(u_1),\dots,f_j(u_m)$ precede $v$, and $f_j(u_s)$ is not one of them. So $f_j(u_s)$ is not among the $b$ smallest of these elements, and by the choice of the subsequence $f_k(u_s)$ is not among the $b$ smallest of the elements $f_k(u_1),\dots,f_k(u_m)$. Since $h(v)$ has colour $d_b$, exactly $b$ of the latter elements precede $h(v)$, and these are the $b$ smallest. Thus $f_k(u_s)$ does not precede $h(v)$. As $h(v)$ and $f_k(u_s)$ have different colours, they are distinct, and therefore $h(v)<_{y_k}f_k(u_s)$.

Finally, let $g\colon x\to y_j$ be a morphism with $g<f_j$, and let $s$ be the least index with $g(u_s)\ne f_j(u_s)$. Then $g(u_s)<_{y_j}f_j(u_s)$. For $t<s$ we have $hg(u_t)=hf_j(u_t)=f_k(u_t)$, and $hg(u_s)<_{y_k}f_k(u_s)$ by the claim. Hence $hg<f_k$.
\end{proof}

\subsection{Proof of Theorem~\ref{thm:criterion}}

\begin{lemma}\label{lem:projective}
Suppose that $\cC$ satisfies the hypothesis of Theorem~\ref{thm:criterion} and that $R$ is left-Noetherian. Then every subrepresentation of $P_x$ is finitely generated.
\end{lemma}

\begin{proof}
Let $M$ be a subrepresentation of $P_x$. For a morphism $f\colon x\to y$, let $I(f)$ be the set of all $a\in R$ for which $M(y)$ contains an element of the form $af+\sum_{g<f}a_gg$, where the sum is over the morphisms $g\colon x\to y$ with $g<f$. Then $I(f)$ is a left ideal of $R$.

We first show that $I(f)\subseteq I(f')$ whenever $f\unlhd f'$. Let $h\colon y\to z$ be a witness and let $a\in I(f)$, so that $M(y)$ contains an element $w=af+\sum_{g<f}a_gg$. Then $M(z)$ contains $P_x(h)(w)=af'+\sum_{g<f}a_ghg$, and $hg<f'$ for every $g<f$. Hence $a\in I(f')$.

Next we claim that there are finitely many morphisms $f_1,\dots,f_n$ with source $x$ and elements $a_j\in I(f_j)$ such that
\[
I(f)=\sum_{j\,:\,f_j\unlhd f}Ra_j
\]
for every morphism $f$ with source $x$. Note that the right-hand side is contained in $I(f)$ for any choice of the $f_j$ and $a_j$, by the previous paragraph. Suppose that the claim is false. Then we can choose morphisms $f_1,f_2,\dots$ and elements $a_j\in I(f_j)$ one after another such that $a_k$ does not belong to the left ideal generated by the elements $a_j$ with $j<k$ and $f_j\unlhd f_k$. By Lemmas~\ref{lem:order} and~\ref{lem:wqo}(i), there are indices $j_1<j_2<\cdots$ with $f_{j_1}\unlhd f_{j_2}\unlhd\cdots$. For every $n$, the element $a_{j_n}$ does not belong to $Ra_{j_1}+\dots+Ra_{j_{n-1}}$. So the left ideals $Ra_{j_1}+\dots+Ra_{j_n}$ form a strictly increasing sequence, contrary to the assumption that $R$ is left-Noetherian. This proves the claim.

Let $f_1,\dots,f_n$ and $a_1,\dots,a_n$ be as in the claim, and write $y_j$ for the target of $f_j$. For each $j$ choose an element $w_j\in M(y_j)$ of the form $a_jf_j+\sum_{g<f_j}a_{j,g}\,g$. Let $M'$ be the subrepresentation of $M$ generated by $w_1,\dots,w_n$. We show that $M'=M$. Suppose not, and choose an object $y$ with $M'(y)\ne M(y)$. Every nonzero element $w$ of $P_x(y)$ has a \emph{leading morphism}, the largest morphism $f\colon x\to y$ which has a nonzero coefficient in $w$. Choose $w\in M(y)\setminus M'(y)$ such that its leading morphism $f$ is as small as possible, and let $a$ be the coefficient of $f$ in $w$. Then $a\in I(f)$, and so $a=\sum_jb_ja_j$ for some elements $b_j\in R$, where the sum is over the indices $j$ with $f_j\unlhd f$. For each of these indices let $h_j\colon y_j\to y$ be a witness for $f_j\unlhd f$, and put
\[
w'=w-\sum_jb_jP_x(h_j)(w_j).
\]
The element $P_x(h_j)(w_j)$ equals $a_jf$ plus a linear combination of morphisms smaller than $f$. So the coefficient of $f$ in $w'$ is $a-\sum_jb_ja_j=0$, and $w'$ is zero or has a leading morphism smaller than $f$. Moreover $w'\in M(y)$ and $w-w'\in M'(y)$, whence $w'\notin M'(y)$. This contradicts the choice of $w$.
\end{proof}

\begin{proof}[Proof of Theorem~\ref{thm:criterion}]
Let $R$ be a left-Noetherian ring. As explained above, we may assume that $\cC$ is small, and the statement for full subcategories follows from the statement for $\cC$, because the hypothesis is inherited by full subcategories. Let $N$ be a finitely generated representation of $\cC$, generated by elements $x_j\in N(y_j)$ for $1\le j\le n$. The morphisms $P_{y_j}\to N$ which map the identity of $y_j$ to $x_j$ define a surjective morphism from $P=P_{y_1}\oplus\dots\oplus P_{y_n}$ onto $N$. The preimage of a subrepresentation of $N$ is a subrepresentation of $P$, so it suffices to show that every subrepresentation of $P$ is finitely generated. We use induction on $n$. The case $n=1$ is Lemma~\ref{lem:projective}. Let $n\ge2$, let $M$ be a subrepresentation of $P$, and let $\pi\colon P\to P_{y_n}$ be the projection. By Lemma~\ref{lem:projective}, the image $\pi(M)$ is generated by finitely many elements, which we lift to elements of $M$. By induction, the kernel of the restriction of $\pi$ to $M$, which is a subrepresentation of $P_{y_1}\oplus\dots\oplus P_{y_{n-1}}$, is finitely generated. The lifted elements together with generators of this kernel generate $M$.
\end{proof}

\subsection{Simplicial complexes}

Many representations of interest are given by the homology of simplicial complexes on the underlying sets. Let $\mathcal S$ be a full subcategory of $\cC$. Suppose that every object $y$ of $\mathcal S$ is assigned a simplicial complex $\Delta(y)$ whose vertices are elements of $V(y)$, and that $h(\sigma)$ is a face of $\Delta(z)$ for every morphism $h\colon y\to z$ of $\mathcal S$ and every face $\sigma$ of $\Delta(y)$. As the maps $V(h)$ are injective, every such morphism induces a chain map between the simplicial chain complexes of $\Delta(y)$ and $\Delta(z)$. Hence $y\mapsto H_i(\Delta(y);R)$ is a representation of $\mathcal S$ over $R$ for every $i\ge0$.

\begin{corollary}\label{cor:complexes}
Suppose that $\cC$ satisfies the hypothesis of Theorem~\ref{thm:criterion}, and let $\mathcal S$ and $\Delta$ be as above. Let $i\ge0$ be an integer.
\begin{enumerate}[label=\textup{(\roman*)}]
\item For every left-Noetherian ring $R$, the representation $y\mapsto H_i(\Delta(y);R)$ of $\mathcal S$ is finitely generated.
\item There is a positive integer $b$ which annihilates the torsion subgroup of $H_i(\Delta(y);\mathbb Z)$ for every object $y$ of $\mathcal S$.
\end{enumerate}
\end{corollary}

\begin{proof}
By Theorem~\ref{thm:criterion}, $\Rep_R(\mathcal S)$ is locally Noetherian for every left-Noetherian ring $R$. Both assertions are invariant under replacing $\mathcal S$ by a skeleton, so, as in Section~\ref{subsec:order}, we may assume that $\mathcal S$ is small. Let $k\ge1$ and write $[k]=\{1,\dots,k\}$. For an object $y$ of $\mathcal S$ let $F_k(y)$ be the free left $R$-module on the set of injections $[k]\to V(y)$, and for a morphism $h\colon y\to z$ let $F_k(h)$ map an injection $\iota$ to $h\iota$. This defines a representation $F_k$ of $\mathcal S$.

We show that $F_k$ is finitely generated. With an injection $\iota\colon[k]\to V(y)$ we associate the colouring $c_\iota$ of $V(y)$ with colours $0,1,\dots,k$ in which $\iota(s)$ has colour $s$ for $1\le s\le k$ and all other elements have colour $0$. For injections $\iota\colon[k]\to V(y)$ and $\kappa\colon[k]\to V(z)$, a morphism $h\colon y\to z$ is colour-preserving from $(y,c_\iota)$ to $(z,c_\kappa)$ if and only if $h\iota=\kappa$. So, by the hypothesis, the pairs $(y,\iota)$ are well-quasi-ordered by the relation in which $(y,\iota)\le(z,\kappa)$ if $h\iota=\kappa$ for some morphism $h\colon y\to z$ of $\mathcal S$. By Lemma~\ref{lem:wqo}(ii), there are finitely many pairs $(y_1,\iota_1),\dots,(y_n,\iota_n)$ such that every pair $(z,\kappa)$ satisfies $(y_j,\iota_j)\le(z,\kappa)$ for some $j$, that is, $\kappa=F_k(h)(\iota_j)$ for some $h\colon y_j\to z$. Hence $\iota_1,\dots,\iota_n$ generate $F_k$.

Let $C_i(y)$ denote the module of oriented simplicial $i$-chains of $\Delta(y)$ with coefficients in $R$. It is generated by the symbols $[v_0,\dots,v_i]$, where $v_0,\dots,v_i$ are the distinct vertices of an $i$-dimensional face in some order, subject to the relations that a permutation of the vertices multiplies the symbol by its sign. Let $F'_{i+1}$ be the subrepresentation of $F_{i+1}$ spanned by the injections whose image is a face. Mapping such an injection $\iota$ to $[\iota(1),\dots,\iota(i+1)]$ defines a surjective morphism of representations from $F'_{i+1}$ onto $C_i$. As $F_{i+1}$ is finitely generated and $\Rep_R(\mathcal S)$ is locally Noetherian, $F'_{i+1}$ is finitely generated, and hence so is $C_i$. The $i$-cycles form a subrepresentation of $C_i$, which is finitely generated for the same reason, and $y\mapsto H_i(\Delta(y);R)$ is a quotient of it. This proves~(i).

For~(ii), let $R=\mathbb Z$. Homomorphisms map torsion elements to torsion elements, so the torsion subgroups $T(y)$ of the groups $H_i(\Delta(y);\mathbb Z)$ form a subrepresentation $T$ of a finitely generated representation. Hence $T$ is generated by finitely many elements $t_1,\dots,t_n$. Let $b$ be a positive integer with $bt_j=0$ for every $j$. Every element of $T(y)$ is an integral combination of images of the elements $t_j$, and is therefore annihilated by $b$.
\end{proof}

\section{A categorical graph minor theorem for edge maps}\label{sec:minors}

In this section we prove Theorem~\ref{thm:edge}. We first consider the injections of edge sets induced by minor models of arbitrary graphs, and then compare them with the injections $\phi^*$ induced by minor morphisms.

\subsection{Edge injections induced by models}

\begin{definition}\label{def:model}
Let $G$ and $H$ be graphs. A \emph{model} of $G$ in $H$ is a pair $(\beta,f)$ with the following properties.
\begin{enumerate}[label=\textup{(\roman*)}]
\item $\beta$ assigns to every vertex $v$ of $G$ a connected subgraph $\beta(v)$ of $H$, called a \emph{branch set}, and $\beta(u)$ and $\beta(v)$ have no common vertex for distinct $u,v\in V(G)$.
\item $f\colon E(G)\to E(H)$ is an injection.
\item For every edge $e$ of $G$ with ends $u$ and $v$, one end of $f(e)$ belongs to $\beta(u)$ and the other end belongs to $\beta(v)$, and $f(e)$ is not an edge of $\beta(u)$ or of $\beta(v)$.
\end{enumerate}
An injection $f\colon E(G)\to E(H)$ is \emph{induced by a model} if $(\beta,f)$ is a model of $G$ in $H$ for some $\beta$.
\end{definition}

In~(iii) the vertices $u$ and $v$ coincide if $e$ is a loop. In this case both ends of $f(e)$ belong to $\beta(v)$, and $f(e)$ may be a loop or an edge with two distinct ends. Since distinct branch sets are disjoint, an edge $f(e)$ as in~(iii) is not an edge of any branch set.

\begin{lemma}\label{lem:category}
Let $G$, $H$ and $K$ be graphs.
\begin{enumerate}[label=\textup{(\roman*)}]
\item The identity of $E(G)$ is induced by a model of $G$ in $G$.
\item If $f\colon E(G)\to E(H)$ and $g\colon E(H)\to E(K)$ are induced by models, then so is $gf$.
\end{enumerate}
\end{lemma}

\begin{proof}
For~(i), let $\beta(v)$ be the subgraph with the single vertex $v$ and no edges.

For~(ii), let $(\beta,f)$ be a model of $G$ in $H$ and $(\gamma,g)$ a model of $H$ in $K$. For a vertex $v$ of $G$, let $\alpha(v)$ be the subgraph of $K$ which is the union of the branch sets $\gamma(w)$ for the vertices $w$ of $\beta(v)$, together with the edges $g(e')$ for the edges $e'$ of $\beta(v)$. Each edge $g(e')$ joins the branch sets of the ends of $e'$. Since $\beta(v)$ and the subgraphs $\gamma(w)$ are connected, $\alpha(v)$ is connected. For distinct vertices $u,v$ of $G$ the subgraphs $\beta(u)$ and $\beta(v)$ have no common vertex, and so $\alpha(u)$ and $\alpha(v)$ have no common vertex. Let $e$ be an edge of $G$ with ends $u$ and $v$. The ends of $f(e)$ are a vertex $w_1$ of $\beta(u)$ and a vertex $w_2$ of $\beta(v)$, and so one end of $gf(e)$ belongs to $\gamma(w_1)\subseteq\alpha(u)$ and the other to $\gamma(w_2)\subseteq\alpha(v)$. The edge $gf(e)$ is not an edge of any $\gamma(w)$. It is not of the form $g(e')$ for an edge $e'$ of a branch set $\beta(w)$ either, because $g$ is injective and $f(e)$ is not an edge of any branch set of $\beta$. Hence $gf(e)$ is not an edge of $\alpha(u)$ or of $\alpha(v)$, and $(\alpha,gf)$ is a model of $G$ in $K$.
\end{proof}

By Lemma~\ref{lem:category}, there is a category $\cE$ whose objects are the graphs and whose morphisms from $G$ to $H$ are the injections $E(G)\to E(H)$ induced by models of $G$ in $H$. It is essentially small, and the functor which maps a graph $G$ to $E(G)$ and a morphism to itself is a faithful functor from $\cE$ to the category of finite sets and injections. We denote by $\Ec$ the full subcategory of $\cE$ on the nonempty connected graphs.

The hypothesis of Theorem~\ref{thm:criterion} for $\cE$ is supplied by the following theorem of Robertson and Seymour. In this statement directed graphs are finite and may have loops and parallel edges.

\begin{theorem}[{\cite[1.7]{RS23}}]\label{thm:RS-minor}
Let $\Omega$ be a well-quasi-ordered set. For every $j\ge1$ let $G_j$ be a directed graph and let $\varphi_j$ be a map from $V(G_j)\cup E(G_j)$ to $\Omega$. Then there are indices $j<k$ and a map $\eta$ with domain $V(G_j)\cup E(G_j)$ such that the following hold.
\begin{enumerate}[label=\textup{(\roman*)}]
\item For every vertex $v$ of $G_j$, $\eta(v)$ is a connected subgraph of $G_k$ which contains a vertex $w$ with $\varphi_j(v)\le\varphi_k(w)$, and $\eta(u)$ and $\eta(v)$ are disjoint for distinct vertices $u,v$.
\item For every edge $e$ of $G_j$, $\eta(e)$ is an edge of $G_k$ with $\varphi_j(e)\le\varphi_k(\eta(e))$, it is a loop if and only if $e$ is a loop, and $\eta(e)\ne\eta(e')$ for distinct edges $e,e'$.
\item If an edge $e$ of $G_j$ has head $u$ and tail $v$, then $\eta(e)$ has its head in $\eta(u)$ and its tail in $\eta(v)$, and $\eta(e)$ is not an edge of $\eta(u)$ or of $\eta(v)$.
\end{enumerate}
\end{theorem}

\begin{proposition}\label{prop:E}
For every finite set $D$, the $D$-coloured objects of $\cE$ are well-quasi-ordered by colour-preserving morphisms. Consequently, $\Rep_R(\cE)$ and $\Rep_R(\mathcal S)$ are locally Noetherian for every left-Noetherian ring $R$ and every full subcategory $\mathcal S$ of $\cE$.
\end{proposition}

\begin{proof}
Let $(G_1,c_1),(G_2,c_2),\dots$ be an infinite sequence of graphs $G_j$ with colourings $c_j\colon E(G_j)\to D$. Orient every edge of every $G_j$ arbitrarily. Let $\Omega$ be the set $D\cup\{*\}$, where $*$ is an element not in $D$, ordered by equality. As $\Omega$ is finite, it is a wqo. Define $\varphi_j(v)=*$ for every vertex $v$ of $G_j$ and $\varphi_j(e)=c_j(e)$ for every edge $e$ of $G_j$. Let $j<k$ and $\eta$ be as in Theorem~\ref{thm:RS-minor}, let $\beta$ be the restriction of $\eta$ to $V(G_j)$ and $f$ the restriction of $\eta$ to $E(G_j)$. Forgetting the orientations, $(\beta,f)$ is a model of $G_j$ in $G_k$ in the sense of Definition~\ref{def:model}, and $c_k(f(e))=c_j(e)$ for every edge $e$ of $G_j$. So $f$ is a colour-preserving morphism of $\cE$ from $(G_j,c_j)$ to $(G_k,c_k)$. Note that Theorem~\ref{thm:RS-minor} represents loops by loops, which is more than Definition~\ref{def:model} requires. The second assertion follows from Theorem~\ref{thm:criterion}.
\end{proof}

\subsection{Minor morphisms and edge injections}

Recall from Section~\ref{sec:trees} that every minor morphism $\phi$ from $H$ to $G$ induces an injection $\phi^*\colon E(G)\to E(H)$.

\begin{lemma}\label{lem:full}
Let $G$ and $H$ be nonempty connected graphs.
\begin{enumerate}[label=\textup{(\roman*)},beginpenalty=10000]
\item For every minor morphism $\phi$ from $H$ to $G$, the injection $\phi^*$ is a morphism from $G$ to $H$ in $\Ec$. Moreover, $(\phi\psi)^*=\psi^*\phi^*$ for composable minor morphisms $\phi$ and $\psi$. Thus $\phi\mapsto\phi^*$ defines a functor $\Phi\colon\Gop\to\Ec$ which is the identity on objects.
\item The functor $\Phi$ is full: for every morphism $f\colon G\to H$ of $\Ec$ there is a minor morphism $\phi$ from $H$ to $G$ with $\phi^*=f$.
\end{enumerate}
\end{lemma}

\begin{proof}
(i) For a vertex $v$ of $G$ let $\beta(v)$ be the tree formed by the vertices and edges of $H$ which are mapped to $v$, as in~(M6). These trees are pairwise disjoint. Let $e$ be an edge of $G$ with ends $u$ and $v$. By~(M4), one end of $\phi^*(e)$ is mapped to $u$ and the other end is mapped to $v$, so the ends of $\phi^*(e)$ belong to $\beta(u)$ and $\beta(v)$. The edge $\phi^*(e)$ is not contracted, so it is not an edge of a branch set. Hence $(\beta,\phi^*)$ is a model of $G$ in $H$. Let $\psi$ be a minor morphism from $K$ to $H$, and let $e$ be an edge of $G$. The arrows of the edge $(\phi\psi)^*(e)$ of $K$ are mapped by $\phi\psi$ to the arrows of $e$. As $\phi$ maps vertices to vertices and $\star$ to $\star$, the map $\psi$ sends them to arrows of $H$, which $\phi$ maps to the arrows of $e$. So $\psi$ maps the edge $(\phi\psi)^*(e)$ to the edge $\phi^*(e)$, that is, $(\phi\psi)^*(e)=\psi^*(\phi^*(e))$.

(ii) Let $(\beta,f)$ be a model of $G$ in $H$, where $G$ and $H$ are nonempty and connected. We first enlarge the branch sets so that they cover $V(H)$. As long as some vertex of $H$ lies in no branch set, there is an edge $xy$ of $H$ such that $x$ lies in a branch set $\beta(v)$ and $y$ lies in none, because $H$ is connected and $G$ has a vertex. We add $y$ and the edge $xy$ to $\beta(v)$. The edge $xy$ is not in the image of $f$, since both ends of every edge $f(e)$ lie in branch sets from the start. So $(\beta,f)$ remains a model. When every vertex of $H$ lies in a branch set, we choose a spanning tree $S_v$ of $\beta(v)$ for every vertex $v$ of $G$. No edge of $S_v$ is in the image of $f$.

Define $\phi$ as follows. Let $\phi(\star)=\star$. For every vertex $v$ of $G$, map the vertices of $\beta(v)$ and both arrows of every edge of $S_v$ to $v$. For every edge $e$ of $G$, map the two arrows of $f(e)$ to the two arrows of $e$ in such a way that heads are mapped to heads and tails to tails; this is possible by Definition~\ref{def:model}(iii), and if $e$ is a loop then either of the two possibilities may be chosen. Map all remaining arrows of $H$ to $\star$. Conditions (M1), (M2), (M4) and~(M5) hold by construction, (M3) holds because $f$ is injective, and~(M6) holds because the vertices and edges mapped to $v$ are those of $S_v$. So $\phi$ is a minor morphism from $H$ to $G$, and $\phi^*=f$.
\end{proof}

By Lemma~\ref{lem:full}, for every representation $\overline M$ of $\Ec$ the composition $\overline M\Phi$ is a $\Gop$-module which depends on a minor morphism $\phi$ only through $\phi^*$. The proof of Theorem~\ref{thm:edge} shows that all $\Gop$-modules with this property arise in this way.

\begin{proof}[Proof of Theorem~\ref{thm:edge}]
By Lemma~\ref{lem:full}, $\Phi$ is the identity on objects and maps the set of minor morphisms from $H$ to $G$ onto the set of morphisms from $G$ to $H$ in $\Ec$. Hence, by the assumption on $M$, there is a unique representation $\overline M$ of $\Ec$ with $M=\overline M\Phi$. The modules $M$ and $\overline M$ have the same values on objects and the same sets of induced homomorphisms. So they have the same submodules, and a family of elements generates $M$ if and only if it generates $\overline M$. The assertion now follows from Proposition~\ref{prop:E}, applied to the full subcategory $\Ec$ of $\cE$.
\end{proof}

\begin{remark}\label{rem:gmt}
Proposition~\ref{prop:E} implies the graph minor theorem. Indeed, let $G_1,G_2,\dots$ be an infinite sequence of graphs, and let $R$ be a field. Let $M$ be the representation of $\cE$ with $M(G)=R$ for every graph $G$, on which every morphism acts as the identity of $R$. Let $G_0$ be the graph with no vertices. The empty injection is a morphism from $G_0$ to every graph, so $M$ is generated by the element $1$ of $M(G_0)$. For $n\ge1$ and a graph $H$, let $N_n(H)=R$ if there is a morphism of $\cE$ from $G_j$ to $H$ for some $j\le n$, and let $N_n(H)=0$ otherwise. Then $N_1\subseteq N_2\subseteq\cdots$ are subrepresentations of $M$, and so is their union $N$. By Proposition~\ref{prop:E}, $N$ is finitely generated, and hence $N=N_n$ for some $n$. As $N(G_{n+1})=R$, there is a morphism of $\cE$ from $G_j$ to $G_{n+1}$ for some $j\le n$, which means that $G_j$ is isomorphic to a minor of $G_{n+1}$. This is not a new proof of the graph minor theorem, since Proposition~\ref{prop:E} rests on Theorem~\ref{thm:RS-minor}.
\end{remark}

\begin{remark}\label{rem:labels}
Theorem~\ref{thm:RS-minor}, applied to the wheels with one label on the spokes and another label on the rim edges, yields integers $k<m$ and a model of $W_k$ in $W_m$ in which spokes are represented by spokes and rim edges by rim edges. Let $f$ be the corresponding injection of edge sets. By Lemma~\ref{lem:full}, there are minor morphisms $\phi$ from $W_m$ to $W_k$ with $\phi^*=f$, but each of them contracts a rim edge of $W_m$. Indeed, let $S_k$ be the set of spokes of $W_k$. If $\phi$ contracted no rim edge, then the spanning tree of $W_m$ with edge set $K(\sigma_k\phi)=K(\phi)\cup\phi^*(S_k)$ would consist of spokes, and so $\sigma_k\phi=\sigma_m$, contrary to Lemma~\ref{lem:wheels}. This is why the labelled graph minor theorem gives the hypothesis of Theorem~\ref{thm:criterion} for $\cE$, but no well-quasi-order on the minor morphisms with target $\bullet$ in which $\tau\le\tau'$ means that $\tau'=\tau\psi$ for some minor morphism $\psi$.
\end{remark}

\section{Matching complexes}\label{sec:matching}

\subsection{Proof of Theorem~\ref{thm:matching}}

\begin{lemma}\label{lem:functorial}
Let $f\colon G\to H$ be a morphism of $\cE$ and let $d\ge1$. If $F$ is a face of $M_d(G)$, then $f(F)$ is a face of $M_d(H)$.
\end{lemma}

\begin{proof}
Let $(\beta,f)$ be a model of $G$ in $H$, and let $w$ be a vertex of $H$. If $w$ belongs to no branch set, then $w$ is not an end of any edge in the image of $f$, and its degree in $f(F)$ is zero. Otherwise $w$ belongs to $\beta(v)$ for exactly one vertex $v$ of $G$. Let $e\in F$, with ends $u_1$ and $u_2$, and let $w_1\in\beta(u_1)$ and $w_2\in\beta(u_2)$ be the corresponding ends of $f(e)$. If $w_1=w$ then $u_1=v$, and if $w_2=w$ then $u_2=v$, because distinct branch sets are disjoint. So the number of ends of $f(e)$ at $w$ is at most the number of ends of $e$ at $v$. Summing over $e\in F$, we see that the degree of $w$ in $f(F)$ is at most the degree of $v$ in $F$, which is at most $d$.
\end{proof}

By Lemma~\ref{lem:functorial}, the complexes $M_d(G)$ satisfy the assumptions of Corollary~\ref{cor:complexes} for the category $\cE$ and the edge-set functor. In particular, $G\mapsto H_i(M_d(G);R)$ is a representation of $\cE$.

\begin{proposition}\label{prop:matching-E}
Let $i\ge0$ be an integer.
\begin{enumerate}[label=\textup{(\roman*)}]
\item For every integer $d\ge1$ and every left-Noetherian ring $R$, the representation $G\mapsto H_i(M_d(G);R)$ of $\cE$ is finitely generated, and so is its restriction to $\Ec$.
\item There is a positive integer $b_i$ such that $b_i$ annihilates the torsion subgroup of $H_i(M_d(G);\mathbb Z)$ for every graph $G$ and every integer $d\ge1$.
\end{enumerate}
\end{proposition}

\begin{proof}
Part~(i) follows from Proposition~\ref{prop:E} and Corollary~\ref{cor:complexes}(i), applied to the full subcategories $\cE$ and $\Ec$ of $\cE$. For~(ii), Corollary~\ref{cor:complexes}(ii) gives, for every $d\ge1$, a positive integer $b_{i,d}$ which annihilates the torsion subgroup of $H_i(M_d(G);\mathbb Z)$ for every graph $G$. Suppose that $d\ge2(i+2)$. In a set of at most $i+2$ edges every vertex has degree at most $2(i+2)$, so every such set is a face of $M_d(G)$. Hence $M_d(G)$ and the full simplex on $E(G)$ have the same faces of dimension at most $i+1$, and therefore the same homology in degree $i$. The homology of a simplex is free, so $H_i(M_d(G);\mathbb Z)$ has no torsion. Thus we may take for $b_i$ the least common multiple of the integers $b_{i,d}$ with $1\le d<2(i+2)$.
\end{proof}

Let $\phi$ be a minor morphism from $H$ to $G$. By Lemmas~\ref{lem:full}(i) and~\ref{lem:functorial}, the injection $\phi^*$ maps the faces of $M_d(G)$ to faces of $M_d(H)$. The homomorphisms of the $\Gop$-modules $G\mapsto H_i(M_d(G);R)$ considered in~\cite{MR} are those induced by the injections $\phi^*$.

\begin{proof}[Proof of Theorem~\ref{thm:matching}]
The $\Gop$-module $G\mapsto H_i(M_d(G);R)$ is the composition of $\Phi$ with the restriction to $\Ec$ of the representation in Proposition~\ref{prop:matching-E}(i). This restriction is finitely generated, and by Lemma~\ref{lem:full}(ii) every morphism of $\Ec$ is of the form $\phi^*$. Hence the same finitely many elements generate the $\Gop$-module, which proves~(i). Part~(ii) is Proposition~\ref{prop:matching-E}(ii).
\end{proof}

\subsection{Other graph complexes}

The argument applies to other complexes of sets of edges. For a graph $G$ and a set $F\subseteq E(G)$, let $G[F]$ denote the subgraph of $G$ formed by the edges in $F$ and their ends. \emph{Splitting} a vertex $v$ of a graph means replacing $v$ by one or more new vertices and distributing the ends of edges at $v$ among the new vertices, so that each new vertex receives at least one of them. Let $\mathcal P$ be a class of graphs without isolated vertices which is closed under isomorphism and contains the graph with no vertices. Suppose that $\mathcal P$ is closed under deleting an edge together with any ends that become isolated, and under splitting a vertex. For a graph $G$, let $K_{\mathcal P}(G)$ be the simplicial complex whose faces are the sets $F\subseteq E(G)$ with $G[F]\in\mathcal P$.

\begin{corollary}\label{cor:properties}
Let $\mathcal P$ be as above and let $i\ge0$ be an integer.
\begin{enumerate}[label=\textup{(\roman*)}]
\item For every left-Noetherian ring $R$, the homology groups $H_i(K_{\mathcal P}(G);R)$ form a finitely generated representation of $\cE$, and their restriction to connected graphs is a finitely generated $\Gop$-module.
\item There is a positive integer $c_i$, depending only on $i$, which annihilates the torsion subgroup of $H_i(K_{\mathcal P}(G);\mathbb Z)$ for every graph $G$ and every class $\mathcal P$ as above.
\end{enumerate}
\end{corollary}

\begin{proof}
Let $f\colon G\to H$ be a morphism of $\cE$, induced by a model $(\beta,f)$, and let $F$ be a face of $K_{\mathcal P}(G)$. The graph $H[f(F)]$ is obtained from $G[F]$ by replacing every edge $e$ by $f(e)$. An end of $e$ at a vertex $v$ becomes an end of $f(e)$ at a vertex of $\beta(v)$, and distinct branch sets are disjoint. So $H[f(F)]$ is obtained from $G[F]$ by splitting vertices, and $f(F)$ is a face of $K_{\mathcal P}(H)$. Thus Corollary~\ref{cor:complexes} applies, and~(i) follows as in the proofs of Proposition~\ref{prop:matching-E} and Theorem~\ref{thm:matching}. Corollary~\ref{cor:complexes} also gives a positive integer $c_{i,\mathcal P}$ which annihilates the torsion subgroup of $H_i(K_{\mathcal P}(G);\mathbb Z)$ for every graph $G$.

Let $\mathcal P_i$ be the class of the members of $\mathcal P$ with at most $i+2$ edges. It satisfies the same assumptions as $\mathcal P$. The complexes $K_{\mathcal P}(G)$ and $K_{\mathcal P_i}(G)$ have the same faces of dimension at most $i+1$, and hence the same homology in degree $i$. Up to isomorphism there are only finitely many graphs without isolated vertices with at most $i+2$ edges, so only finitely many classes occur as $\mathcal P_i$. The least common multiple of the integers $c_{i,\mathcal P_i}$ for these classes has the property required in~(ii).
\end{proof}

\begin{remark}\label{rem:properties}
The graphs without isolated vertices in which every vertex has degree at most $d$ form a class $\mathcal P$ as above, for which $K_{\mathcal P}(G)=M_d(G)$. Other examples are the forests and the bipartite graphs without isolated vertices. If a loop is counted only once in the degree of its vertex, one obtains a variant of the $d$-matching complex which is also of the form $K_{\mathcal P}$, so that Corollary~\ref{cor:properties} applies to it as well. For graphs of bounded first Betti number, Miyata and Ramos~\cite[Proposition~4.5]{MR} proved finite generation for the complexes of monotone graph properties which are compatible with minor morphisms.
\end{remark}

\section{Immersions}\label{sec:immersions}

In this section graphs are loopless. Let $G$ and $H$ be loopless graphs. An \emph{immersion} of $G$ in $H$ is a map $\alpha$ which assigns to every vertex $v$ of $G$ a vertex $\alpha(v)$ of $H$ and to every edge $e$ of $G$ with ends $u$ and $v$ a path $\alpha(e)$ of $H$ with ends $\alpha(u)$ and $\alpha(v)$, such that $\alpha(u)\ne\alpha(v)$ for distinct vertices $u,v$, and the paths $\alpha(e)$ and $\alpha(e')$ have no common edge for distinct edges $e,e'$. The restriction of $\alpha$ to $V(G)$ is the \emph{vertex map} of the immersion. Immersions can be composed: given an immersion of $G$ in $H$ and an immersion of $H$ in $K$, replace every edge of the path representing an edge of $G$ by the path of $K$ representing it, and choose a path inside the resulting walk. The vertex map of the composition is the composition of the vertex maps.

\begin{theorem}[{\cite[1.5]{RS23}}]\label{thm:RS-immersion}
Let $\Omega$ be a well-quasi-ordered set. For every $j\ge1$ let $G_j$ be a loopless graph and let $\varphi_j$ be a map from $V(G_j)$ to $\Omega$. Then there are indices $j<k$ and an immersion $\alpha$ of $G_j$ in $G_k$ such that $\varphi_j(v)\le\varphi_k(\alpha(v))$ for every vertex $v$ of $G_j$.
\end{theorem}

Let $Q$ be a well-quasi-ordered set. A \emph{$Q$-labelled graph} is a loopless graph $G$ together with a map $\lambda_G\colon V(G)\to Q$. Let $\mathcal J_Q$ be the category whose objects are the $Q$-labelled graphs and whose morphisms from $G$ to $H$ are the vertex maps of the immersions $\alpha$ of $G$ in $H$ with $\lambda_G(v)\le\lambda_H(\alpha(v))$ for every vertex $v$ of $G$. The vertex sets define a faithful functor from $\mathcal J_Q$ to the category of finite sets and injections.

\begin{corollary}\label{cor:immersion}
For every well-quasi-ordered set $Q$ and every left-Noetherian ring $R$, the categories $\Rep_R(\mathcal J_Q)$ and $\Rep_R(\mathcal S)$ are locally Noetherian for every full subcategory $\mathcal S$ of $\mathcal J_Q$.
\end{corollary}

\begin{proof}
Let $D$ be a finite set. Order $Q\times D$ by $(q,a)\le(q',a')$ if $q\le q'$ and $a=a'$. This is a wqo, since an infinite sequence in $Q\times D$ has an infinite subsequence with constant second coordinate. A $D$-coloured object of $\mathcal J_Q$ is a loopless graph with vertex labels in $Q\times D$, and a colour-preserving morphism is the vertex map of an immersion which respects these labels. By Theorem~\ref{thm:RS-immersion}, the $D$-coloured objects of $\mathcal J_Q$ are well-quasi-ordered by colour-preserving morphisms, and the assertion follows from Theorem~\ref{thm:criterion}.
\end{proof}

Corollary~\ref{cor:complexes} now applies to simplicial complexes on vertex sets whose faces are preserved by the vertex maps of immersions. For example, call a set $S$ of vertices of a loopless graph $G$ a face if the complete graph on $S$ has an immersion in $G$ whose vertex map is the identity of $S$. Subsets of faces are faces, and the vertex map of an immersion of $G$ in $H$ maps faces to faces. If $G$ is the path with three vertices, then every set of two vertices is a face, but the set of all three vertices is not, since it would require three paths with at least one edge each and no common edges, while $G$ has only two edges. So the complex of $G$ is the boundary of a triangle, and its first homology group is $\mathbb Z$.

\section{Concluding remarks}\label{sec:conclusion}

\subsection*{Bounds on torsion}

The proof of Theorem~\ref{thm:matching} relies on the well-quasi-ordering theorem of Robertson and Seymour and is not effective. It gives neither an explicit set of generators in part~(i) nor a bound on $b_i$ in part~(ii). As noted in the introduction, every admissible $b_1$ is divisible by $3$.

\begin{question}\label{q:torsion}
What is the least integer $b_1$ with the property in Theorem~\ref{thm:matching}(ii)? More generally, find an explicit upper bound for the least possible value of $b_i$.
\end{question}

\subsection*{Minor-closed classes}

For a class $\mathcal K$ of graphs which is closed under taking minors, let $\cG_{\mathcal K}$ be the full subcategory of $\cG$ on the connected graphs in $\mathcal K$.

\begin{question}\label{q:classes}
For which minor-closed classes $\mathcal K$ of graphs is every submodule of a finitely generated $\cG_{\mathcal K}^{\mathrm{op}}$-module over a Noetherian ring finitely generated?
\end{question}

By~\cite[Theorem~1.1]{MR} this is the case for the class of graphs all of whose components have first Betti number at most $g$, for every $g$. By Remark~\ref{rem:wheels} it is not the case when $\mathcal K$ contains all wheels, for instance for the class of planar graphs, and by Remark~\ref{rem:brooms} it is not the case for the class of graphs with no minor isomorphic to $K_4$.

\subsection*{Conditional results}

Some finiteness results in the literature are conditional on Conjecture~\ref{conj:cgm}; see for instance~\cite[Proposition~6.5 and Corollary~6.6]{KR}, and~\cite[Theorem~1.3]{KR}, which assumes a second conjecture as well. Theorem~\ref{thm:trees} shows that they need a different justification, and Proposition~\ref{prop:torsion} shows that the torsion of a finitely generated $\Gop$-module need not be bounded. Theorem~\ref{thm:edge} provides such a justification for all $\Gop$-modules which depend on a minor morphism $\phi$ only through the injection $\phi^*$. It would be interesting to know which of the modules arising in these applications have this property, or are quotients of modules that have it.

\section*{Use of AI tools}

During this research and the preparation of the manuscript, the authors used Anthropic's Claude Fable~5.2 and OpenAI's ChatGPT (6.1 Sol) in order to discuss the proof strategies. These discussions contributed to the development of the finite-colouring criterion and its applications, and to the identification of the examples in Section~\ref{sec:trees}: the wheels used in the proof of Theorem~\ref{thm:trees}, the module of Proposition~\ref{prop:torsion}, and the graphs of Remark~\ref{rem:brooms}. The tools were also used to check proofs and references, to run computer checks of small cases, and to improve the exposition. The authors have checked all statements, proofs and references, have reviewed and edited the text, and take full responsibility for the content of this article.

\end{document}